\documentclass[11pt]{amsart}

\usepackage{latexsym}
\usepackage{amssymb}

\newcommand{\excise}[1]{}%{$\star$\textsc{#1}$\star$}

\newtheorem{thm}{Theorem}%[section]
\newtheorem{lemma}[thm]{Lemma}

\newtheorem{cor}[thm]{Corollary}
\newtheorem{prop}[thm]{Proposition}
\newtheorem{conj}[thm]{Conjecture}

\newtheorem{Example}[thm]{Example}
\newtheorem{Remark}[thm]{Remark}
\newtheorem{Alg}[thm]{Algorithm}
\newtheorem{Defn}[thm]{Definition}

\newenvironment{example}{\begin{Example}\rm}
                {\mbox{}~\hfill$\square$\end{Example}}
\newenvironment{remark}{\begin{Remark}\rm}
                {\mbox{}~\hfill$\square$\end{Remark}}

\newenvironment{defn}{\begin{Defn}\rm}{\end{Defn}}

\newenvironment{defnlabeled}[1]{\begin{Defn}[#1]\rm}{\end{Defn}}

\newenvironment{rcgraph}{\begin{trivlist}\item\centering\footnotesize$}
                        {$\end{trivlist}}

\def\hln{\\[-.2ex]\hline}

\newenvironment{numbered}%
        {\begin{list}
                {\noindent\makebox[0mm][r]{\arabic{enumi}.}}
                {\leftmargin=5.5ex \usecounter{enumi}}
        }
        {\end{list}}

        {\begin{list}
                {\noindent\makebox[0mm][r]{(\roman{enumi})}}
                {\leftmargin=5.5ex \usecounter{enumi}}
        }
        {\end{list}}

\newcounter{separated}

\def\bem#1{\textbf{#1}}

\def\<{\langle}
\def\>{\rangle}
\def\0{{\mathbf 0}}
\def\1{{\mathbf 1}}

\def\JJ{{\mathcal J}}

\def\LL{{\mathcal L}}

\def\SS{{\mathfrak S}}

\def\ZZ{{\mathbb Z}}

\def\xx{{\mathbf x}}

\def\id{{\rm id}}

\def\th{{\rm th}}

\def\cross{\textrm{`+'} }

\def\start{{\rm start}}

\def\facets{{\rm facets}}
\def\length{{\rm length}}

\def\mitosis{{\rm mitosis\hspace{1pt}}}

\def\rc{\mathcal{RP}}
\def\dom{\backslash}

\def\adots{{.\hspace{1pt}\raisebox{2pt}{.}\hspace{1pt}\raisebox{4pt}{.}}}
\def\minus{\smallsetminus}

\def\nothing{\varnothing}
\def\bigcupdot{\makebox[0pt][l]{$\hspace{1.05ex}\cdot$}\textstyle\bigcup}
\def\toplinedots{\hfill\raisebox{-5.1pt}[0pt][0pt]{\ $\ldots$\ }\hfill}

\def\rcn#1{\mathcal{RP}_{\!#1}}

\font\co=lcircle10
\def\petit#1{{\scriptstyle #1}}

\def\dt{\hspace{.17ex}\cdot\hspace{.17ex}}
\def\jr{\smash{\raise2pt\hbox{\co \rlap{\rlap{\char'005} \char'007}}
               \raise6pt\hbox{\rlap{\vrule height6.5pt}}
               \raise2pt\hbox{\rlap{\hskip4pt \vrule height0.4pt depth0pt
                width7.7pt}}}}
\def\je{\smash{\raise2pt\hbox{\co \rlap{\rlap{\char'005}
                \phantom{\char'007}}}\raise6pt\hbox{\rlap{\vrule height6pt}}}}
\def\+{\smash{\lower2pt\hbox{\rlap{\vrule height14pt}}
                \raise2pt\hbox{\rlap{\hskip-3pt \vrule height.4pt depth0pt
                width14.7pt}}}}
\def\perm#1#2{\hbox{\rlap{$\petit {#1}_{\scriptscriptstyle #2}$}}%
                \phantom{\petit 1}}
\def\phperm{\phantom{\perm w3}}

\def\reflection#1#2#3{\raisebox{#1}{${\scriptstyle#2}_{\scriptscriptstyle#3}$}}

\def\textcross{\ \smash{\lower4pt\hbox{\rlap{\hskip4.15pt\vrule height14pt}}
                \raise2.8pt\hbox{\rlap{\hskip-3pt \vrule height.4pt depth0pt
                width14.7pt}}}\hskip12.7pt}
\def\textelbow{\ \hskip.1pt\smash{\raise2.8pt%
                \hbox{\co \hskip 4.15pt\rlap{\rlap{\char'005} \char'007}
                \lower6.8pt\rlap{\vrule height3.5pt}
                \raise3.6pt\rlap{\vrule height3.5pt}}
                \raise2.8pt\hbox{%
                  \rlap{\hskip-7.15pt \vrule height.4pt depth0pt width3.5pt}%
                  \rlap{\hskip4.05pt \vrule height.4pt depth0pt width3.5pt}}}
                \hskip8.7pt}
\def\plel{$\begin{tinyrc}{
  \begin{array}{@{}|@{\,}c@{\,}|@{}}
      \hline \scriptscriptstyle  +
    \\\hline \cdot
    \\\hline
  \end{array}
  }\end{tinyrc}$\ }
\def\elel{$\begin{tinyrc}{
  \begin{array}{@{}|@{\,}c@{\,}|@{}}
      \hline \,\cdot\,
    \\\hline \cdot
    \\\hline
  \end{array}
  }\end{tinyrc}$\ }
\def\elpl{$\begin{tinyrc}{
  \begin{array}{@{}|@{\,}c@{\,}|@{}}
      \hline \cdot
    \\\hline \scriptscriptstyle   +  
    \\\hline
  \end{array}
  }\end{tinyrc}$\ }
\def\plpl{$\begin{tinyrc}{
  \begin{array}{@{}|@{\,}c@{\,}|@{}}
      \hline \scriptscriptstyle   +  
    \\\hline \scriptscriptstyle   +  
    \\\hline
  \end{array}
  }\end{tinyrc}$\ }

\begin{document}%%%%%%%%%%%%%%%%%%%%%%%%%%%%%%%%%%%%%%%%%%%%%%%%%%%%%%%%%
%%%%%%%%%%%%%%%%%%%%%%%%%%%%%%%%%%%%%%%%%%%%%%%%%%%%%%%%%%%%%%%%%%%%%%%%%

\title{Mitosis recursion for coefficients of Schubert polynomials}
\author{Ezra Miller}
\thanks{The author was supported by the Sloan Foundation and NSF}
\address{MSRI\\1000 Centennial Drive\\Berkeley, CA}
\email{emiller@msri.org}
\date{5 December 2002}

\begin{abstract}
\noindent
Mitosis is a rule introduced in \cite{grobGeom} for manipulating subsets
of the $n \times n$ grid.  It provides an algorithm that lists the
reduced pipe dreams (also known as rc-graphs) \cite{FKyangBax,BB} for a
permutation \mbox{$w \in S_n$} by downward induction on weak Bruhat
order, thereby generating the coefficients of Schubert polynomials
\cite{LSpolySchub} inductively.  This note provides a short and purely
combinatorial proof of these properties of mitosis.
%\vskip 1ex
%\noindent
%{{\it AMS Classification:} ; }
\end{abstract}

\maketitle

%%%%%%%%%%%%%%%%%%%%%%%%%%%%%%%%%%%%%%%%%%%%%%%%%%%%%%%%%%%%%%%%%%%%%%%%%
{}%%%%%%%%%%%%%%%%%%%%%%%%%%%%%%%%%%%%%%%%%%%%%%%%%%%%%%%%%%%%%%%%%%%%%%%
%%%%%%%%%%%%%%%%%%%%%%%%%%%%%%%%%%%%%%%%%%%%%%%%%%%%%%%%%%%%%%%%%%%%%%%%%

%%%%%%%%%%%%%%%%%%%%%%%%%%%%%%%%%%%%%%%%%%%%%%%%%%%%%%%%%%%%%%%%%%%%%%%%%
\section{Introduction}%%%%%%%%%%%%%%%%%%%%%%%%%%%%%%%%%%%%%%%%%%%%%%%%%%%
%%%%%%%%%%%%%%%%%%%%%%%%%%%%%%%%%%%%%%%%%%%%%%%%%%%%%%%%%%%%%%%%%%%%%%%%%

It has been a goal for some years, ever since Kohnert made his
conjecture in \cite{Kohnert}, to find inductive combinatorial rules on
diagrams in the $n \times n$ grid that yield the coefficients of
Schubert polynomials \cite{LSpolySchub}, when counted properly.  The
mitosis rule was offered in \cite{grobGeom} as a solution to this
problem, but the proof was long, and involved some notions that strayed
rather far from the elementary combinatorics of permutations.  The
purpose of this note is to bring mitosis entirely into the realm of
combinatorics, by giving a short combinatorial proof of the fact
(Theorem~\ref{thm:mitosis}) that mitosis lists reduced pipe dreams (also
known as rc-graphs) \cite{FKyangBax,BB} recursively by induction on weak
order in~$S_n$, starting from the unique reduced pipe dream for the long
permutation~$w_0$.

More precisely, the proof here of Theorem~\ref{thm:mitosis}, and the
resulting diagrammatic recursion for the coefficients of Schubert
polynomials in Corollary~\ref{cor:mitosis}, rests only on the formula of
Billey, Jockusch, and Stanley (Theorem~\ref{thm:BJS}), the
characterization of Schubert polynomials by divided differences
(Definition~\ref{defn:schub}), and elementary combinatorial properties of
reduced pipe dreams (Lemmas~\ref{lemma:rc}, \ref{lemma:chute},
and~\ref{lemma:intron} plus Proposition~\ref{prop:sym}).

Mitosis was originally conceived in \cite{grobGeom} as a residual
operation derived from more complicated combinatorial isobaric divided
differences (Demazure operators) on standard monomials for certain
determinantal ideals defined in the context of Schubert varieties in
flag manifolds.  As such, it served as a geometrically motivated
improvement on Kohnert's rule \cite{Kohnert, NoteSchubPoly, Winkel,
Winkel02}, its advantages being the short combinatorial proof here and
consistency with double Schubert polynomials as in~\cite{grobGeom}.
Mitosis is closely related to the construction of Schubert polynomials
in terms of chains in Bruhat order in~\cite{LenSot02}.  Other
combinatorial algorithms producing Schubert polynomials include the
chute and ladder moves on reduced pipe dreams \cite{BB}, a different
combinatorial divided difference operator on reduced pipe dreams
\cite{LenUnified}, and an earlier construction of Bergeron~\cite{Ber92}.

The plan of the paper is as follows.  In the next two sections we review
the definition of the set $\rc(w)$ of reduced pipe dreams for a
permutation $w \in S_n$, the BJS formula, and the mitosis algorithm on
pipe dreams (subsets of the $n \times n$ grid).  Section~\ref{sec:intron}
provides an involution on~$\rc(w)$ that is crucial for the proof of the
main theorem and corollary in Section~\ref{sec:thm}.  The final section,
which concerns the mitosis poset and is logically independent of the
other sections, reviews for the reader's convenience two definitions and
a conjecture from \cite[Section~2.2]{grobGeom}, because of their
relevance in this combinatorial setting.

%%%%%%%%%%%%%%%%%%%%%%%%%%%%%%%%%%%%%%%%%%%%%%%%%%%%%%%%%%%%%%%%%%%%%%%%%
\section{Pipe dreams}%%%%%%%%%%%%%%%%%%%%%%%%%%%%%%%%%%%%%%%%%%%%%%%%%%%%
%%%%%%%%%%%%%%%%%%%%%%%%%%%%%%%%%%%%%%%%%%%%%%%%%%%%%%%%%%%%%%%%%%%%%%%%%

\label{sec:pipe}

Consider a square grid $\ZZ_{> 0} \times \ZZ_{> 0}$ extending infinitely
south and east, with the box in row~$i$ and column~$j$ labeled $(i,j)$,
as in an $\infty \times \infty$ matrix.  If each box in the grid is
covered with a square tile containing either $\textcross$ or
$\textelbow$, then one can think of the tiled grid as a network of pipes.

\begin{defn} \label{defn:pipe}
A \bem{pipe dream} is a finite subset of $\ZZ_{> 0} \times \ZZ_{> 0}$,
identified as the set of crosses in a tiling by \bem{crosses}
$\textcross$ and \bem{elbow joints} $\textelbow$.  A pipe dream is
\bem{reduced} if each pair of pipes crosses at most once.  The set
$\rc(w)$ of reduced pipe dreams for the permutation $w \in S_n$ is the
set of reduced pipe dreams~$D$ such that the pipe entering row~$i$ exits
from column $w(i)$.%
\end{defn}

Although we always draw crossing tiles as some sort of cross (either
\cross or `$\textcross$', the former with the square tile boundary and
the latter without), we often leave the elbow tiles blank or denote them
by dots, to make the diagrams less cluttered.  Viewing $n$ as fixed, we
shall be interested in pipe dreams contained in the pipe dream~$D_0$
that has crosses in the triangular region strictly above the main
antidiagonal (in spots $(i,j)$ with $i+j \leq n$) and elbow joints
elsewhere in the square grid $[n] \times [n]$ of size~$n$.  Note that
$D_0$ is the unique reduced pipe dream for the \bem{long permutation}
$w_0 = n \ldots 3 2 1$ in~$S_n$.

\begin{example} \label{ex:pipe}
The pipe dream~$D$ in Fig.~\ref{fig:pipe} for $n = 8$ is a reduced pipe
dream for the permutation $w = 13865742 \in S_8$.
\begin{figure}[ht]
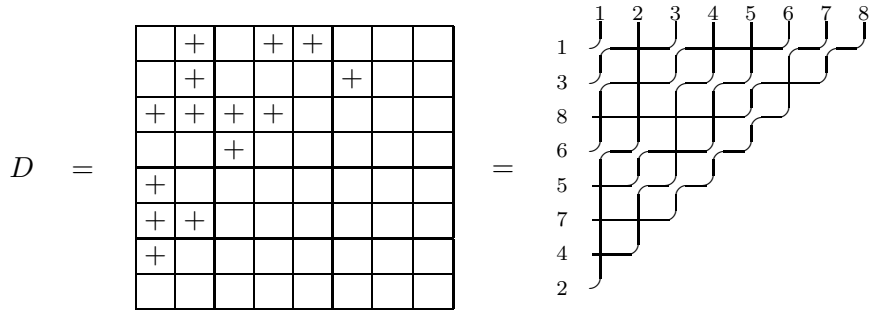

$$
\begin{array}{@{}c@{}}\\[-6ex]
\begin{array}{c}\\\\\\\\\\\raisebox{2ex}{$D\quad=\quad$}\\\\\\\\\end{array}
\begin{array}{|c|c|c|c|c|c|c|c|}
\multicolumn{8}{c}{}
\\\hline     &\!+\!&     &\!+\!&\!+\!&     &\ \, &\ \,
\\\hline     &\!+\!&     &     &     &\!+\!&     &
\\\hline\!+\!&\!+\!&\!+\!&\!+\!&     &     &     &
\\\hline     &     &\!+\!&     &     &     &     &
\\\hline\!+\!&     &     &     &     &     &     &
\\\hline\!+\!&\!+\!&     &     &     &     &     &
\\\hline\!+\!&     &     &     &     &     &     &
\\\hline     &     &     &     &     &     &     &
\\\hline
\end{array}
\ \quad \raisebox{-1.5ex}{$=$} \quad
\begin{array}{ccccccccc}
      &\petit1&\petit 2&\petit 3&\petit 4&\petit 5&\petit 6&\petit 7&\petit 8\\
\petit1&  \jr &   \+   &   \jr  &   \+   &   \+   &   \jr  &   \jr  &  \je\\
\petit3&  \jr &   \+   &   \jr  &   \jr  &   \jr  &   \+   &   \je  &\\
\petit8&  \+  &   \+   &   \+   &   \+   &   \jr  &   \je  &        &\\
\petit6&  \jr &   \jr  &   \+   &   \jr  &   \je  &        &        &\\
\petit5&  \+  &   \jr  &   \jr  &   \je  &        &        &        &\\
\petit7&  \+  &   \+   &   \je  &        &        &        &        &\\
\petit4&  \+  &   \je  &        &        &        &        &        &\\
\petit2&  \je &        &        &        &        &        &        &
\end{array}
\\[-2ex]
\end{array}
$$
\caption{A reduced pipe dream for $w = 13865742 \in S_8$}
\label{fig:pipe}
\end{figure}
For clarity, we omit the square tile boundaries as well as the wavy
``sea'' of elbows $\textelbow$ below the main antidiagonal in the right
pipe dream.%
% The thinner symbol $w_i$ instead of $w(i)$ makes the column widths come
% out right.%
\end{example}

Since we need a statement of the BJS formula, we recall here the
definition of Schubert polynomials of Lascoux and Sch\"utzenberger via
divided differences.  For notation, $s_i \in S_n$ denotes the
transposition switching $i$ and~\mbox{$i+1$}, and $\length(w)$ denotes
the number of inversions in a permutation~$w$.

\begin{defnlabeled}{\cite{LSpolySchub}} \label{defn:schub}
The \bem{$i^\th$ divided difference operator} $\partial_i$ takes each
polynomial $f \in \ZZ[x_1,\ldots,x_n]$~to
\begin{eqnarray*}
  \partial_i f(x_1,\ldots,x_n) &=& \frac{f(x_1,\ldots,x_n) - f(x_1,
  \ldots, x_{i-1}, x_{i+1}, x_i, x_{i+2}, \ldots, x_n)}{x_i - x_{i+1}}.
\end{eqnarray*}
The \bem{Schubert polynomial} for~$w \in S_n$ is defined by the recursion
\begin{eqnarray*}
  \SS_{ws_i}(x_1,\ldots,x_n) &=& \partial_i \SS_w(x_1,\ldots,x_n)
\end{eqnarray*}
whenever $\length(ws_i) < \length(w)$, and the initial condition
$\SS_{w_0}(x_1,\ldots,x_n) = \prod_{i=1}^n x_i^{n-i}$.%
\end{defnlabeled}

\begin{thm}[\cite{BJS,FSnilCoxeter}] \label{thm:BJS}
$\displaystyle \SS_w(x_1,\ldots,x_n)\ = \sum_{D \in \rc(w)} \xx^D$,\ \
where $\ \displaystyle \xx^D = \prod_{(i,j) \in D} x_i$.
\end{thm}

The next lemma, which will be applied in Section~\ref{sec:thm}, gives a
criterion for when removing a~\cross from a pipe dream $D \in \rc(w)$
leaves a pipe dream in~$\rc(ws_i)$.  Specifically, it concerns the
removal of a cross at $(i,j)$ from configurations that look like
\begin{rcgraph}
\begin{array}{lccccccccc}
&\perm 1{}&\perm{}{}&\perm{}{}&\perm{}{}&
        \perm{\hskip-11.6pt \textstyle \cdots}{}&\perm{}{}&\perm{}{}&\perm j\
\\[3pt]
\petit{i}  &\+  &  \+  &  \+  &  \+  &  \+  &  \+  &  \+  &  \+ &\toplinedots
\\[-1pt]
\petit{i+1}&\+  &  \+  &  \+  &  \+  &  \+  &  \+  &  \+  &  \jr
\end{array}
\begin{array}{c}\\ = \\\end{array}\
\begin{array}{l|c|c|c|c|c|c|c|c|c}
\multicolumn{1}{c}{}
        &\multicolumn{1}{c}{\petit 1}
                &\multicolumn{6}{c}{\cdots}
                        &\multicolumn{1}{c}{\petit j}
                                \\[3pt]\cline{2-10}
\petit{i}  &\!+\!&\!+\!&\!+\!&\!+\!&\!+\!&\!+\!&\!+\!&\!+\!&\toplinedots
                                                            \\\cline{2-9}
\petit{i+1}&\!+\!&\!+\!&\!+\!&\!+\!&\!+\!&\!+\!&\!+\!&\cdot&\\\cline{2-10}
\end{array}
\end{rcgraph}
at the left end of rows~$i$ and~$i+1$ in~$D$.

\begin{lemma} \label{lemma:rc}
Let $D \in \rc(w)$ and $j$ be a fixed column index with $(i+1,j) \not\in
D$, but $(i,p) \in D$ for all $p \leq j$, and $(i+1,p) \in D$ for all $p
< j$.  Then $\length(ws_i) < \length(w)$, and if $D' = D \minus (i,j)$
then $D' \in \rc(ws_i)$.
\end{lemma}
\begin{proof}
Removing $(i,j)$ only switches the exit points of the two pipes starting
in rows $i$ and $i+1$, so the pipe starting in row $k$ of $D'$ exits out
of column $ws_i(k)$ for each $k$.  No pair of pipes can cross twice in
$D'$ because there are $\length(ws_i)$ crossings.%
\end{proof}

%%%%%%%%%%%%%%%%%%%%%%%%%%%%%%%%%%%%%%%%%%%%%%%%%%%%%%%%%%%%%%%%%%%%%%%%%
\section{Mitosis algorithm}%%%%%%%%%%%%%%%%%%%%%%%%%%%%%%%%%%%%%%%%%%%%%%
%%%%%%%%%%%%%%%%%%%%%%%%%%%%%%%%%%%%%%%%%%%%%%%%%%%%%%%%%%%%%%%%%%%%%%%%%

\label{sec:alg}

Given a pipe dream in $[n] \times [n]$, define
\begin{eqnarray} \label{eq:pipestart}
  \start_i(D) &=& \hbox{column index of leftmost elbow in row } i\\
\nonumber      &=& \min(\{j \mid (i,j) \not\in D\} \cup \{n+1\}),
\end{eqnarray}
so the $i^\th$ row of $D$ is filled solidly with crosses in the region to
the left of $\start_i(D)$.  Let
\begin{eqnarray*}
  \JJ_i(D) &=& \{\hbox{columns $j$ strictly to the left of } \start_i(D)
  \mid (i+1,j) \hbox{ has no cross in } D\}.
\end{eqnarray*}
For $p \in \JJ_i(D)$, construct the \bem{offspring} $D_p$ as follows.
First delete the cross at $(i,p)$ from~$D$.  Then take all crosses in
row~$i$ of $\JJ_i(D)$ that are to the left of column~$p$, and move each
one down to the empty box below it in row~$i+1$.

\begin{defn} \label{defn:mitosis}
The $i^\th$ \bem{mitosis} operator sends a pipe dream $D$ to
\begin{eqnarray*}
  \mitosis_i(D) &=& \{D_p \mid p \in \JJ_i(D)\}.
\end{eqnarray*}
Write $\mitosis_i({\mathcal P}) = \bigcup_{D \in {\mathcal P}}
\mitosis_i(D)$ whenever ${\mathcal P}$ is a set of pipe dreams.%
\end{defn}

Observe that all of the action takes place in rows~$i$ and~$i+1$, and
$\mitosis_i(D)$ is an empty set whenever $\JJ_i(D)$ is empty.

\begin{example} \label{ex:mitosis}
The pipe dream $D$ at left is the reduced pipe dream for $w = 13865742$
from Example~\ref{ex:pipe}:
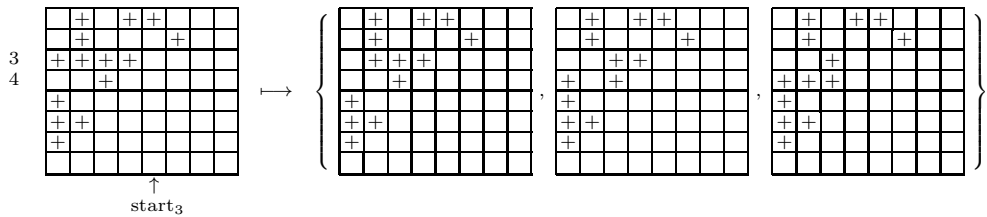
$$
\begin{tinyrc}{
\begin{array}{@{}cc@{}}{}\\\\3&\\4\\\\\\\\\\\\\\\end{array}
\begin{array}{@{}|@{\,}c@{\,}|@{\,}c@{\,}|@{\,}c@{\,}|@{\,}c@{\,}
		 |@{\,}c@{\,}|@{\,}c@{\,}|@{\,}c@{\,}|@{\,}c@{\,}|@{}}
\hline     &  +  &     &  +  &  +  &     &     &
\hln       &  +  &     &     &     &  +  &     &\phantom{+}
\hln    +  &  +  &  +  &  +  &     &     &     &
\hln       &     &  +  &     &     &     &     &
\hln    +  &     &     &     &     &     &     &
\hln    +  &  +  &     &     &     &     &     &
\hln    +  &     &     &     &     &     &     &
\hln       &     &     &     &     &     &\phantom{+}&
\\\hline\multicolumn{4}{@{}c@{}}{}&\multicolumn{1}{@{}c@{}}{\!\!\uparrow}
\\\multicolumn{4}{c}{}&\multicolumn{1}{@{}c@{}}{\makebox[0pt]{$\start_3$}}
\end{array}
\begin{array}{@{\quad}c@{\ }}
  \longmapsto\\\mbox{}\\\mbox{}
\end{array}
\begin{array}{c}
  \left\{\ 
  \begin{array}{@{}|@{\,}c@{\,}|@{\,}c@{\,}|@{\,}c@{\,}|@{\,}c@{\,}
  		 |@{\,}c@{\,}|@{\,}c@{\,}|@{\,}c@{\,}|@{\,}c@{\,}|@{}}
  \hline     &  +  &     &  +  &  +  &     &     &    
  \hln       &  +  &     &     &     &  +  &     &\phantom{+}
  \hln       &  +  &  +  &  +  &     &     &     &
  \hln       &     &  +  &     &     &     &     &
  \hln    +  &     &     &     &     &     &     &
  \hln    +  &  +  &     &     &     &     &     &
  \hln    +  &     &     &     &     &     &     &
  \hln       &     &     &     &     &     &\phantom{+}&
  \\\hline
  \end{array}
  \ ,\ 
  \begin{array}{@{}|@{\,}c@{\,}|@{\,}c@{\,}|@{\,}c@{\,}|@{\,}c@{\,}
  		 |@{\,}c@{\,}|@{\,}c@{\,}|@{\,}c@{\,}|@{\,}c@{\,}|@{}}
  \hline     &  +  &     &  +  &  +  &     &     &    
  \hln       &  +  &     &     &     &  +  &     &\phantom{+}
  \hln       &     &  +  &  +  &     &     &     &
  \hln    +  &     &  +  &     &     &     &     &
  \hln    +  &     &     &     &     &     &     &
  \hln    +  &  +  &     &     &     &     &     &
  \hln    +  &     &     &     &     &     &     &
  \hln       &     &     &     &     &     &\phantom{+}&
  \\\hline
  \end{array}
  \ ,\ 
  \begin{array}{@{}|@{\,}c@{\,}|@{\,}c@{\,}|@{\,}c@{\,}|@{\,}c@{\,}
  		 |@{\,}c@{\,}|@{\,}c@{\,}|@{\,}c@{\,}|@{\,}c@{\,}|@{}}
  \hline     &  +  &     &  +  &  +  &     &     &    
  \hln       &  +  &     &     &     &  +  &     &\phantom{+}
  \hln       &     &  +  &     &     &     &     &
  \hln    +  &  +  &  +  &     &     &     &     &
  \hln    +  &     &     &     &     &     &     &
  \hln    +  &  +  &     &     &     &     &     &
  \hln    +  &     &     &     &     &     &     &
  \hln       &     &     &     &     &     &\phantom{+}&
  \\\hline
  \end{array}
  \ \right\}
\\
  \begin{array}{c}
  \mbox{}\\
  \mbox{}\\
  \end{array}
\end{array}
}\end{tinyrc}
$$
The set of three pipe dreams on the right is obtained by applying
$\mitosis_3$, since $\JJ_3(D)$ consists of columns $1$, $2$, and $4$.
The offspring are ordered as in Proposition~\ref{prop:offspring}, below.%
\end{example}

In Proposition~\ref{prop:offspring} we shall present another, more
sequential way of writing down the mitosis offspring of a pipe dream.  It
uses a device invented by Bergeron and Billey.

\begin{defnlabeled}{\cite{BB}} \label{defn:chute}
A \bem{chutable rectangle} is a connected $2 \times k$ rectangle $C$
inside a pipe dream $D$ such that $k \geq 2$ and all but the following 3
locations in $C$ are crosses: the northwest, southwest, and southeast
corners.  Applying a \bem{chute move} to $D$ is accomplished by placing a
\cross in the southwest corner of a chutable rectangle $C$ and removing
the \cross from the northeast corner of the same~$C$.%
\end{defnlabeled}

Heuristically, a chute move therefore looks like:
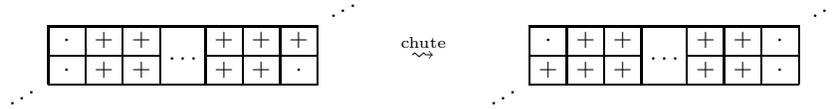
\begin{rcgraph}
\begin{array}{@{}c@{}}\\[-5ex]
%                1 2 3 4    5    6 7 8 9
\begin{array}{@{}r|c|c|c|@{}c@{}|c|c|c|l@{}}
  \multicolumn{7}{c}{}&\multicolumn{1}{c}{
                                          \phantom{\!+\!}}&
  \multicolumn{1}{c}{\begin{array}{@{}c@{}}\\\adots\end{array}}
  \\\cline{2-8}
           &\cdot&\!+\!&\!+\!& \toplinedots &\!+\!&\!+\!&\!+\!
  \\\cline{2-4}\cline{6-8}
           &\cdot&\!+\!&\!+\!&              &\!+\!&\!+\!&\cdot
  \\\cline{2-8}
  \multicolumn{1}{c}{\begin{array}{@{}c@{}}\adots\\ \\ \end{array}}&
  \multicolumn{1}{c}{\phantom{\!+\!}}
\end{array}
\quad\stackrel{\rm chute}\rightsquigarrow\quad
%
%                1 2 3 4    5    6 7 8 9
\begin{array}{@{}r|c|c|c|@{}c@{}|c|c|c|l@{}}
  \multicolumn{7}{c}{}&\multicolumn{1}{c}{
                                          \phantom{\!+\!}}&
  \multicolumn{1}{c}{\begin{array}{@{}c@{}}\\\adots\end{array}}
  \\\cline{2-8}
           &\cdot&\!+\!&\!+\!& \toplinedots &\!+\!&\!+\!&\cdot
  \\\cline{2-4}\cline{6-8}
           &\!+\!&\!+\!&\!+\!&              &\!+\!&\!+\!&\cdot
  \\\cline{2-8}
  \multicolumn{1}{c}{\begin{array}{@{}c@{}}\adots\\ \\ \end{array}}&
  \multicolumn{1}{c}{\phantom{\!+\!}}
\end{array}
\\[-3ex]
\end{array}
\end{rcgraph}

The following basic fact about chute moves was discovered by Bergeron and
Billey~\cite{BB}.

\begin{lemma} \label{lemma:chute}
The set $\rc(w)$ of reduced pipe dreams for~$w$ is closed under chute
moves.
\end{lemma}
\begin{proof}
If two pipe intersect at the \cross in the northeast corner of a chutable
rectangle~$C$, then chuting that \cross only changes the crossing point
of the two pipes to the southwest corner of~$C$.  No other pipes are
affected.%
\end{proof}

\begin{prop} \label{prop:offspring}
Let $D$ be a pipe dream, and suppose $j$ is the smallest column index
such that $(i+1,j) \not\in D$ and $(i,p) \in D$ for all $p \leq j$.  Then
$D_p \in \mitosis_i(D)$ is obtained from~$D$~by
\begin{numbered}
\item
removing $(i,j)$, and then
\item
performing chute moves from row~$i$ to row~$i+1$, each one as far left as
possible, so that $(i,p)$ is the last {\rm \cross$\!$} removed.
\end{numbered}
\end{prop}
\begin{proof}
Immediate from Definitions~\ref{defn:mitosis} and~\ref{defn:chute}.
\end{proof}

%%%%%%%%%%%%%%%%%%%%%%%%%%%%%%%%%%%%%%%%%%%%%%%%%%%%%%%%%%%%%%%%%%%%%%%%%
\section{Intron mutation}%%%%%%%%%%%%%%%%%%%%%%%%%%%%%%%%%%%%%%%%%%%%%%%%
%%%%%%%%%%%%%%%%%%%%%%%%%%%%%%%%%%%%%%%%%%%%%%%%%%%%%%%%%%%%%%%%%%%%%%%%%

\label{sec:intron}

\begin{defn} \label{defn:intron}
Let $D$ be a pipe dream and $i$ a fixed row index.  Order the boxes in
rows $i$ and~$i+1$ of~$D$ as in the following diagram:
\begin{rcgraph}
\begin{array}{l|c|c|c|c|c}
\multicolumn{6}{@{}c@{}}{}\\[-3ex]
\multicolumn{1}{l}{\mbox{}}
  &\multicolumn{1}{c}{\petit{1}}
    &\multicolumn{1}{c}{\petit{2}}
      &\multicolumn{1}{c}{\petit{3}}
        &\multicolumn{1}{c}{\petit{4}}
          &\multicolumn{1}{c}{\cdots}
                                        \\[.2ex]\cline{2-6}
\petit{i}  &  1  &  3  &  5  &  7  &\toplinedots        \\\cline{2-5}
\petit{i+1}&  2  &  4  &  6  &  8  &                    \\\cline{2-6}
\multicolumn{6}{@{}c@{}}{}\\[-1ex]
\end{array}
\end{rcgraph}
An \bem{intron}%
	\footnote{For the origin of this term, see
	\cite[Section~3.5]{grobGeom}.}
in these two adjacent rows is a $2 \times k$ rectangle $C$ such that
\begin{numbered}
\item
the first and last boxes in $C$ (the northwest and southeast corners) are
elbows; and

\item
no elbow in $C$ is strictly northeast or strictly southwest of another
elbow (so due north, due south, due east, or due west are all okay).
\end{numbered}\setcounter{separated}{\value{enumi}}
Ignoring all \plpl columns in rows $i$ and~\mbox{$i+1$}, an intron is
just a sequence of \elpl columns in rows $i$ and~\mbox{$i+1$}, followed
by a sequence of \plel columns, possibly with one \elel column in
between.  Columns with two crosses \plpl can be ignored for the purpose
of proofs in what follows.

If an intron $C$ satisfies the following extra condition, then $C$ is
called a \bem{maximal} intron:
\begin{numbered}\setcounter{enumi}{\value{separated}}
\item
the elbow with largest index before $C$ (if there is one) resides in
row~$i+1$, and the elbow with smallest index after $C$ (if there is one)
resides in row~$i$.
\end{numbered}
\end{defn}

\begin{lemma} \label{lemma:intron}
For an intron $C$ in a reduced pipe dream, a unique intron $\tau(C)$
satisfies
\begin{numbered}
\item
the sets of columns with exactly two crosses are the same in $C$ and
$\tau(C)$, and

\item
the number $c_i$ of crosses in row~$i$ of\/ $C$ equals the number of
crosses in row~$i+1$ of~$\tau(C)$, and conversely.
\end{numbered}
The involution $\tau$, called \bem{intron mutation}, is always
accomplished by a sequence of chute moves or inverse chute moves (because
$C$ is part of a\/ {\em reduced} pipe dream).
\end{lemma}
\begin{proof}
First assume $c_i > c_{i+1}$ and work by induction on $c = c_i -
c_{i+1}$.  If $c = 0$ then $\tau(C) = C$ and the lemma is obvious.  If $c
> 0$ then consider the leftmost \plel column.  Moving to~the left from
this column there must be a column not equal to \plpl\!, since the
northwest entry of~$C$ is an elbow.  The rightmost such column must be
\elel\!, because its row~$i$ entry is an elbow (by construction) and its
row~\mbox{$i+1$} entry cannot be a cross (for then the pipes crossing
there would also cross in the \plel column).  This means we can chute the
\cross in \plel into the \elel column, and proceed by induction.

Flip the argument $180^\circ$ if $c_i < c_{i+1}$, so the chute move
becomes an inverse chute.%
\end{proof}

For example, here is an intron mutation accomplished by chuting the
crosses in columns $4$, $6$, and then $7$ of row~$i$; the zigzag shapes
formed by the dots in these introns are typical.
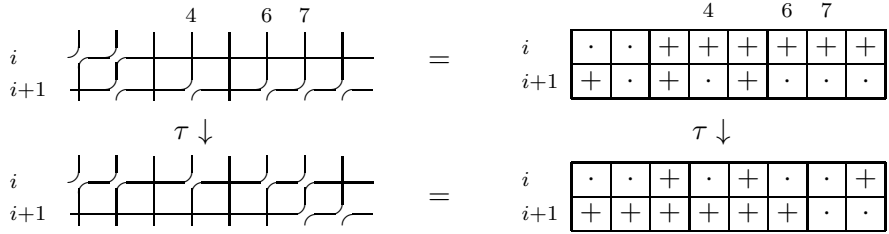
\begin{eqnarray*}
\begin{array}{c@{\ \qquad}c@{\qquad}c}
\begin{array}{lcccccccc}
\multicolumn{8}{c}{}\\[-4ex]
        &\phperm&\phperm&\phperm&\petit4&\phperm&\petit6&\petit7&\phperm\\[3pt]
\petit{i} &  \jr &  \jr &  \+   &  \+   &  \+   &  \+   &  \+   &  \+
\\[-1pt]
\petit{i+1}& \+  &  \jr &  \+   &  \jr  &  \+   &  \jr  &  \jr  &  \jr
\end{array}
&
=
&
\begin{array}{l|c|c|c|c|c|c|c|c|}
\multicolumn{9}{@{}c@{}}{}\\[-5ex]
\multicolumn{4}{l}{\mbox{}}
  &\multicolumn{1}{c}{\petit{4}}
    &\multicolumn{1}{c}{\petit{}}
      &\multicolumn{1}{c}{\petit{6}}
        &\multicolumn{1}{c}{\petit{7}}
          &\multicolumn{1}{c}{}
                                        		\\[.2ex]\cline{2-9}
\petit{i}  &\cdot& \dt &\!+\!&\!+\!&\!+\!&\!+\!&\!+\!&\!+\!
                                                                \\\cline{2-9}
\petit{i+1}&\!+\!&\cdot&\!+\!&\cdot&\!+\!&\cdot&\cdot&\cdot     \\\cline{2-9}
\end{array}
\\
\begin{array}{c}\\[-2ex]\quad\tau \downarrow\\[1ex]\end{array}
&&
\begin{array}{c}\\[-2ex]\quad\tau \downarrow\\[1ex]\end{array}
\\
\begin{array}{lcccccccc}
\multicolumn{8}{c}{}\\[-6ex]
        &\phperm&\phperm&\phperm&\phperm&\phperm&\phperm&\phperm&\phperm\\[3pt]
\petit{i} &  \jr &  \jr &  \+   &  \jr  &  \+   &  \jr  &  \jr  &  \+
\\[-1pt]
\petit{i+1}& \+  &  \+  &  \+   &  \+   &  \+   &  \+   &  \jr  &  \jr
\end{array}
&
=
&
\begin{array}{l|c|c|c|c|c|c|c|c|}
                                                                \cline{2-9}
\petit{i}  &\cdot&\cdot&\!+\!&\cdot&\!+\!&\cdot& \dt &\!+\!
                                                                \\\cline{2-9}
\petit{i+1}&\!+\!&\!+\!&\!+\!&\!+\!&\!+\!&\!+\!&\cdot&\cdot     \\\cline{2-9}
\end{array}
\end{array}
\end{eqnarray*}

\begin{prop} \label{prop:sym}
For each $i$ there is an involution $\tau_i : \rc(w) \to \rc(w)$ such
that $\tau_i^2 = 1$, and for all $D \in \rc(w)$:
\begin{numbered}
\item
$\tau_iD$ agrees with $D$ outside rows $i$ and $i+1$.

\item
$\start_i(\tau_iD) = \start_i(D)$, and $\tau_iD$ agrees with $D$ strictly
west of this column.

\item
$\ell^i_i(\tau_iD) = \ell^i_{i+1}(D)$,
\end{numbered}
where $\ell^i_r(-)$ is the number of crosses in row $r$ that are east of
or in column $\start_i(-)$.
\end{prop}
\begin{proof}
Let $D \in \rc(w)$.  Consider the union of all columns in rows $i$ and
$i+1$ of $D$ that are east of or coincide with column $\start_i(D)$.
Since the first and last boxes in this region (numbered as in
Definition~\ref{defn:intron}) are elbows, this region breaks uniquely
into a disjoint union of $2 \times k$ rectangles, each of which is either
a maximal intron or completely filled with crosses.  Indeed, this follows
from (\ref{eq:pipestart}) and Definition~\ref{defn:intron}.  Applying
intron mutation to each maximal intron therein leaves a pipe dream that
breaks up uniquely into maximal introns and solid crosses in the same
way.  Therefore the lemma comes down to verifying that intron mutation
preserves the property of being in~$\rc(w)$, which comes from
Lemmas~\ref{lemma:chute} and~\ref{lemma:intron}.%
% More accurately, mutation of a single intron in rows $i$ and $i+1$
% of~$D$ yields a pipe dream in $\rc(w)$.  To verify this statement, one
% need only check the routes of pipes intersecting the intron, and this
% is straightforward: the numbers of crosses traversed horizontally by
% each pipe are the same in $C$ and $\tau(C)$; similarly for the
% vertically traversed crosses.%
\end{proof}

\begin{excise}{%
\begin{remark} \label{rk:sym}
The intron mutation illustrated in (\ref{eq:intron}) results from a
sequence of chute moves.  This phenomenon is general, whenever the pipe
dream is reduced: in any single intron of a reduced pipe dream, either
the mutation or its inverse results from a sequence of chute moves.  The
proof is omitted, but rests on the fact that the $2 \times 2$
configuration
$\begin{tinyrc}{
  \begin{array}{@{}|@{\,}c@{\,}|@{\,}c@{\,}|@{}}
      \hline \cdot &  +
    \\\hline   +   &\cdot
    \\\hline
  \end{array}
  }\end{tinyrc}$
is disallowed in reduced pipe dreams, so every intron in an reduced pipe
dream contains a column with two elbows.%
\end{remark}
}\end{excise}%
\begin{remark}
Intron mutation is precisely the involution (coplactic
operation)~$\sigma_i$ defined by Lascoux on words (see the survey article
\cite{LLT}, for example) and extended to reduced pipe dreams
in~\cite{LenUnified}.  However, when all introns in rows $i$
and~\mbox{$i+1$} are strung together, the involution~$\tau_i$ does not
agree with~$\sigma_i$.  In fact, Lascoux's involution is based on
`$r$-pairing', which is also used in the work of Bergeron~\cite{Ber92}
and Lenart~\cite{LenUnified} to define combinatorial versions of divided
difference operators.  Intron mutation is therefore a different mechanism
by which combinatorial divided differences can be defined on reduced pipe
dreams.%
\end{remark}

%%%%%%%%%%%%%%%%%%%%%%%%%%%%%%%%%%%%%%%%%%%%%%%%%%%%%%%%%%%%%%%%%%%%%%%%%
\section{Mitosis theorem}%%%%%%%%%%%%%%%%%%%%%%%%%%%%%%%%%%%%%%%%%%%%%%%%
%%%%%%%%%%%%%%%%%%%%%%%%%%%%%%%%%%%%%%%%%%%%%%%%%%%%%%%%%%%%%%%%%%%%%%%%%

\label{sec:thm}

\begin{thm} \label{thm:mitosis}
If\/ $\length(ws_i) < \length(w)$, then the set $\rc(ws_i)$ of reduced
pipe dreams for~$ws_i$ is the disjoint union $\bigcupdot_{D \in \rc(w)}
\mitosis_i(D)$.  Therefore
\begin{eqnarray} \label{eq:alg}
  \rc(w) &=& \mitosis_{i_k} \cdots \mitosis_{i_1}(D_0)
\end{eqnarray}
if $s_{i_1} \cdots s_{i_k}$ is a reduced expression for~$w_0w$.
\end{thm}
\begin{proof}
Use the description of mitosis in Proposition~\ref{prop:offspring} along
with Lemmas~\ref{lemma:rc} and~\ref{lemma:chute} to conclude that
$\mitosis_i(D) \subseteq \rc(ws_i)$ whenever $D \in \rc(w)$.  It follows
directly from the definitions that $\mitosis_i(D) \cap \mitosis_i(D') =
\nothing$ if $D \neq D'$ are reduced pipe dreams for~$w$.  Thus it
suffices to prove that $\mitosis_i(\rc(w))$ has the same cardinality as
$\rc(ws_i)$.

Fix $D \in \rc(w)$, write $\xx^D = \prod_{(i,j) \in D} x_i$, and let $J =
|\JJ_i(D)|$ be the number of mitosis offspring of~$D$.  The monomial
$\xx^D$ is a product $x_i^J \xx^{D'}$, where $D'$ is the pipe dream ({\em
not}\/ reduced) obtained from~$D$ by erasing the crosses in row $i$ of
$\JJ_i(D)$.  Definition~\ref{defn:mitosis} implies~that
\begin{equation} \label{eq:mitosis}
\sum_{E \in \mitosis_i(D)} \xx^E
\:\ =\:\ \sum_{d=1}^J x_i^{J-d}x_{i+1}^{d-1} \cdot \xx^{D'}
\:\ =\:\ \partial_i(x_i^J) \cdot \xx^{D'}.
%\:\ =\:\ \partial_i(x_i^J \cdot \xx^{D'})
%\:\ =\:\ \partial_i(\xx^D).
\end{equation}
If $\tau_i D = D$, then $\xx^{D'}$ is symmetric in $x_i$ and $x_{i+1}$ by
Proposition~\ref{prop:sym}, so that
$$
\partial_i(x_i^J) \cdot \xx^{D'}
\:\ =\:\ \partial_i(x_i^J \cdot \xx^{D'})
\:\ =\:\ \partial_i(\xx^D)
$$
in this case.  On the other hand, if $\tau_i D \neq D$, then letting
$s_i$ act on polynomials by switching $x_i$ and $x_{i+1}$,
Proposition~\ref{prop:sym} implies that adding the sums in
(\ref{eq:mitosis}) for $D$ and $\tau_i D$ yields
$$
\partial_i(x_i^J) \cdot (\xx^{D'} + s_i \xx^{D'})
\:\ =\:\ \partial_i(x_i^J(\xx^{D'} + s_i\xx^{D'})) 
\:\ =\:\ \partial_i(\xx^D + \xx^{\tau_i D}).
$$
Pairing off the elements of $\rc(w)$ not fixed by~$\tau_i$, we therefore
conclude that
$$
  \sum_{E \in \mitosis_i(\rc(w))} \xx^E \ \:=\:\ \partial_i \Bigl(\sum_{D
  \in \rc(w)} \xx^D \Bigr) \ \:=\:\ \partial_i(\SS_w(\xx)) \ \:=\:\
  \SS_{ws_i}(\xx) \ \:=\: \sum_{E \in \rc(ws_i)} \xx^E
$$
by Theorem~\ref{thm:BJS} and the recursion for $\SS_w(\xx) :=
\SS_w(x_1,\ldots,x_n)$ as in Definition~\ref{defn:schub}.  Plugging in
$1, \ldots, 1$ for $\xx = x_1,\ldots,x_n$ implies that
$|\mitosis_i(\rc(w))| = |\rc(ws_i)|$, as desired.%
\end{proof}

Finally we come to the generation of Schubert coefficients by induction
on weak Bruhat order via mitosis.  For notation, if $v = s_{i_1} \cdots
s_{i_k}$ is a reduced expression, set $\mitosis_v = \mitosis_{i_k} \cdots
\mitosis_{i_1}$.

\begin{cor} \label{cor:mitosis}
For any permutation $w \in S_n$ we have
\begin{eqnarray*}
  \displaystyle \SS_w(x_1,\ldots,x_n) &=& \sum_{D \in {\mitosis_v(D_0)}}
  \xx^D \qquad{\rm for\ } v = w_0w,
\end{eqnarray*}
where $\rc(w_0) = \{D_0\}$, and $\xx^D = \prod_{(i,j) \in D} x_i$ for any
pipe dream~$D$.
\end{cor}
\begin{proof}
Theorem~\ref{thm:mitosis} and Theorem~\ref{thm:BJS}.
\end{proof}

% comment{This corollary works also for double, but it's kinda hard to
% get at combinatorially without resorting to Demazure operators.  You
% see, it's the $\star$elbows$\star$ that get counted in \cite{grobGeom},
% and really we'd like to think of elbows as counting for arbitrary
% positive integers in order to get things right.  What I'm trying to say
% is that the preservation of column sums under intron transposition just
% doesn't work in the purely combinatorial context of reduced pipe dreams
% because a chute move really coalesces two elbows (creating a cross in
% the bottom row) and then ``uncoalesces'' two elbows (erasing a cross
% in the top row).  Thus it is not as easy to work directly with double
% Schubert polynomials in the combinatorial context as it is to work
% with double Grothendieck polynomials}

%%%%%%%%%%%%%%%%%%%%%%%%%%%%%%%%%%%%%%%%%%%%%%%%%%%%%%%%%%%%%%%%%%%%%%%%%
\section{Mitosis poset}%%%%%%%%%%%%%%%%%%%%%%%%%%%%%%%%%%%%%%%%%%%%%%%%%%
%%%%%%%%%%%%%%%%%%%%%%%%%%%%%%%%%%%%%%%%%%%%%%%%%%%%%%%%%%%%%%%%%%%%%%%%%

\label{sec:mitosis}

The next definition generalizes to arbitrary $n$ the poset of pipe
dreams for $n=3$ in Fig.~\ref{fig:hasse}.

\begin{defnlabeled}{{\cite[Definition~2.2.4]{grobGeom}}}
Theorem~\ref{thm:mitosis} defines a partial order, namely
\begin{eqnarray*}
  D' \prec D &\hbox{ if }& D' \in \mitosis_i(D) \hbox{ for some } i,
\end{eqnarray*}
making the reduced pipe dreams for all of~$S_n$ into the \bem{mitosis
poset} $\rcn n = \bigcup_{w \in S_n} \rc(w)$.
\end{defnlabeled}

The poset $\rcn n$, which is ranked by length = cardinality, fibers over
the weak Bruhat order on $S_n$, with the preimage of $w \in S_n$ being
$\rc(w)$.  A reduced expression for $w_0w$ can be thought of as the edge
labels on a decreasing path beginning at $w_0$ and ending at~$w$ in the
weak Bruhat order on~$S_n$.  The preimage in $\rcn n$ of such a path is
a tree having $\rc(w)$ among its leaves (two reduced pipe dreams cannot
share an offspring by the disjointness of the union in
Theorem~\ref{thm:mitosis}).%
\begin{figure}[t]
$$
\begin{array}{@{}r@{}c@{}l@{}}
\\[-4ex]
&
  \hbox{\footnotesize 321}
\\[2pt]
&
  \begin{tinyrc}{
  \begin{array}{@{}|@{\,}c@{\,}|@{\,}c@{\,}|@{\,}c@{\,}|@{}}
      \hline + & + &\phantom{+}
    \\\hline + &   &
    \\\hline   &   &
    \\\hline
  \end{array}
  }\end{tinyrc}
&
\\
  \reflection{5pt} s2\diagup\,
&
&
  \,\diagdown \reflection{5pt} s1
\\
  \begin{array}{@{}c@{}}
    \makebox[0pt][r]{\hbox{\footnotesize 312\ \ }}
    \begin{tinyrc}{
    \begin{array}{@{}|@{\,}c@{\,}|@{\,}c@{\,}|@{\,}c@{\,}|@{}}
        \hline + & + &\phantom{+}
      \\\hline   &   &
      \\\hline   &   &
      \\\hline
    \end{array}
    }\end{tinyrc}
  \\
    \reflection{2pt} s1/ \ \: \dom \reflection{2pt} s1
  \\
    \makebox[0pt][r]{\hbox{\footnotesize 132\ \ }}
    \begin{tinyrc}{
    \begin{array}{|@{\,}c@{\,}|@{\,}c@{\,}|@{\,}c@{\,}|}
        \hline\phantom{+}& + &\phantom{+}
      \\\hline           &   &
      \\\hline           &   &
      \\\hline
    \end{array}
    }\end{tinyrc}
    \:
    \begin{tinyrc}{
    \begin{array}{|@{\,}c@{\,}|@{\,}c@{\,}|@{\,}c@{\,}|}
        \hline   &\phantom{+}&\phantom{+}
      \\\hline + &           &
      \\\hline   &           &
      \\\hline
    \end{array}
    }\end{tinyrc}
  \end{array}
&
&
  \phantom{\,\diagdown}
  \begin{array}{@{\,}c@{}}
    \begin{tinyrc}{
    \begin{array}{@{}|@{\,}c@{\,}|@{\,}c@{\,}|@{\,}c@{\,}|@{}}
        \hline + &\phantom{+}&\phantom{+}
      \\\hline + &           &
      \\\hline   &           &
      \\\hline
    \end{array}
    }\end{tinyrc}
    \makebox[0pt][l]{\hbox{\footnotesize \ \ 231}}
  \\
    |\makebox[0pt][l]{$\reflection{2pt} s2$}
  \\
    \begin{tinyrc}{
    \begin{array}{@{}|@{\,}c@{\,}|@{\,}c@{\,}|@{\,}c@{\,}|@{}}
        \hline + &\phantom{+}&\phantom{+}
      \\\hline   &           &
      \\\hline   &           &
      \\\hline
    \end{array}
    }\end{tinyrc}
    \makebox[0pt][l]{\hbox{\footnotesize \ \ 213}}
  \end{array}
\\
  \reflection{-1pt} s2\diagdown\,
&
&
  \,\diagup \reflection{-1pt} s1
\\
&
  \begin{tinyrc}{
  \begin{array}{@{}|@{\,}c@{\,}|@{\,}c@{\,}|@{\,}c@{\,}|@{}}
      \hline \phantom{+}&\phantom{+}&\phantom{+}
    \\\hline            &           &
    \\\hline            &           &
    \\\hline
  \end{array}
  }\end{tinyrc}
\\[7pt]
&
  \hbox{\footnotesize 123}
\\[-2ex]
\end{array}
$$
\caption{Hasse diagram for~$\rcn 3$} \label{fig:hasse}
\end{figure}%

% the efficiency of the algorithm can vary widely with the reduced
% expression $s_{i_1} \cdots s_{i_k}$ for $w_0w$.  For instance, the set
% $\rc(\id_n)$ for the identity permutation consists of one element,
% $\nothing \subset [n]^2$, even though the repeated mitosis
% in~(\ref{eq:alg}) for $w = \id_n$ can pass through $\rc(w)$ for {\em
% any}\/ permutation $w$.  As a consequence, huge numbers of reduced
% pipe dreams may be killed along the way, without producing any
% offspring. 

\begin{defn} \label{defn:poptotic}
A path decreasing from~$w_0$ to~$w$ in the weak order is \bem{poptotic}
if the leaves of its preimage in $\rcn n$ are precisely $\rc(w)$.%
\end{defn}

In other words, a path is poptotic if every reduced pipe dream lying over
its interior has at least one offspring.  For example, the right hand
path in Fig.~\ref{fig:hasse} from $321$ to $123$ is poptotic because only
one reduced pipe dream appears at each stage, while the left path is
apoptotic%
        \footnote{The word `apoptosis' refers in biology to programmed
        cell death, where some cell in a multicellular organism commits
        suicide for the greater good of the organism.  Thus
        \bem{apoptotic} indicates that some reduced pipe dream dies
        without offspring, while \bem{poptotic} indicates that all pipe
        dreams survive~with~offspring.}
because the first reduced pipe dream for $132$ has no offspring
under~$\mitosis_2$.

\begin{prop} \label{p:poptotic}
Poptotic paths from $w_0$ to~$w$ exist.  In fact, the lexicographically
first reduced expression for $w_0w$ (in which $s_1 > s_2 > \cdots >
s_{n-1}$) corresponds to a poptotic path.
\end{prop}
In particular, the lex first path from $w_0$ to~$\id_n$ passes through
\bem{dominant} permutations, which by definition have exactly one reduced
pipe dream (shaped like a Young diagram).

\begin{proof}
% The lex first reduced expression for $w_0w$ starts with~$w_0$, puts $n$
% and~\mbox{$n-1$} in the same order as in~$w$, then puts $\{n,n-1,n-2\}$
% in the same order as they are in~$w$, and so on.
Number the boxes in the strict upper-left triangle, meaning all locations
$(q,p)$ such that $q+p \leq n$, as follows, where $N = \binom n2$.
\def\mcn#1{\multicolumn{#1}{@{}c@{}}{}} \def\phaN{\phantom{8}}
\def\drop{\begin{array}{@{}c@{}}\vdots\\[-5ex]\end{array}}
\begin{rcgraph}
\begin{array}{|c|c|c|c|c|c|c|}
\hline
\!\!N\!\!&        &\!\!10\!\!& 6  &  3  &  1  &\ \:
\\\cline{1-1}\cline{3-6}
  \phaN  &        &  9  &  5  &  2  &\mcn1&
\\\cline{3-5}
  \drop  & \cdots &  8  &  4  &\mcn2&
\\\cline{3-4}
         &        &  7  &\mcn3&
\\\cline{3-3}
         &\multicolumn{1}{|c@{}}{\raisebox{1ex}{$\!\!\!\!\adots$}}&\mcn4&
\end{array}
\end{rcgraph}
The ordered sequence $(1,2,1,3,2,1,4,3,2,1,\ldots)$ of row indices of
boxes in this upper triangle gives rise to the lex first reduced
expression $s_1s_2s_1s_3s_2s_1s_4s_3s_2s_1\cdots\ $ for the long
word~$w_0 = w_0 \id_n$.  In general, lex first reduced words for
arbitrary~$w_0w$ correspond bijectively to the complements in the
upper-left triangle of so-called \bem{top} reduced pipe dreams \cite{BB},
which are characterized (by definition) as having no \elpl
configurations.  The reduced word corresponding to a top pipe dream is
the ordered subsequence of row indices skipping the crosses atop each
column.

Now suppose that $\length(ws_i) < \length(w)$, and that the lex first
reduced expression for~$w_0ws_i$ ends in~$s_{i'}s_i$.  Under the
bijection above between lex first reduced words and complements of top
reduced pipe dreams, $s_{i'}$ and~$s_i$ correspond to the row indices
$i'$ and~$i$ of boxes numbered~$\alpha'$ and~$\alpha$ satisfying $\alpha'
< \alpha$.  [N.B.\ Either $i = i'-1$, in which case $\alpha' = \alpha-1$,
or else $i > i'$, and $\alpha$ sits just above the main antidiagonal in
some column to the left of~$\alpha'$.]  Therefore, we shall assume by
induction on length that
\setcounter{separated}{\value{equation}}%
\renewcommand{\theequation}{\fnsymbol{equation}} \setcounter{equation}{0}%
\begin{equation}
\hbox{every reduced pipe dream in~$\rc(w)$ has crosses in boxes
$>\alpha'$ and an elbow at~$\alpha'$.}
\end{equation}
\setcounter{equation}{\value{separated}}%
(the case $\length(w_0w) = 1$ is easy).  The goal is to prove that~$(*)$
holds with $ws_i$ in place of~$w$ and $\alpha$ in place of~$\alpha'$.
But first, note that Lemma~\ref{lemma:rc}, which holds with $\alpha$ in
position $(i,j)$ by assumption~$(*)$, says that $\mitosis_i(D)$ is
nonempty for all $D \in \rc(w)$, as required.

More precisely, Lemma~\ref{lemma:rc} says that removing the cross
at~$\alpha$ from each $D \in \rc(w)$ produces a pipe dream
in~$\rc(ws_i)$.  Furthermore, either $\alpha$ lies in the top row or the
box in~$D$ due north of~$\alpha$ is a cross, so it is impossible for
chute moves to end there after deleting the cross from~$\alpha$.
Consequently, Proposition~\ref{prop:offspring} implies that every pipe
dream $D' \in \mitosis_i(D)$ has crosses in boxes marked $>\alpha$, and
an elbow joint in the box marked~$\alpha$.  The proof is complete by
Theorem~\ref{thm:mitosis}.%
\end{proof}

\begin{example} \label{ex:apoptosis}
The three pipe dreams on the right in Example~\ref{ex:mitosis} are all
reduced pipe dreams for $v = 13685742 = w \cdot s_3$, where $w =
13865742$ as in Example~\ref{ex:pipe}.  Setting $i = 4$ and inspecting
the inversions of~$v$, we find that $\length(vs_4) < \length(v)$.  On
the other hand, $\mitosis_4$ kills the first two of the three pipe
dreams, whereas the last has two offspring.  Thus any path from $w_0$
to~$vs_4$ ending with $(\ldots, v, vs_4)$ is necessarily apoptotic.

Note that the lex first reduced expression for $w_0v$, which corresponds
to a poptotic path from $w_0$ to~$v$ by Proposition~\ref{p:poptotic},
equals $s_2s_1s_3s_5s_4s_3s_2s_1s_7s_6s_5s_4s_3s_2s_1$, while the lex
first reduced expression for $w_0vs_4$ equals
$s_2s_1s_3s_2s_5s_4s_3s_2s_1s_7s_6s_5s_4s_3s_2s_1$ (the $s_2$ in the
fourth slot is new).  These correspond to top reduced pipe dreams%
% \begin{verbatim}
% 87654321	
% 87654321	.
% 68754321	. 21
% 65874321	. 21 32 
% 65874321	. 21 32 . 
% 36587421	. 21 32 . 54321
% 36587421	. 21 32 . 54321 . 
% 13658742	. 21 32 . 54321 . 7654321
% 
% 13658742
% 
% 87654321	.
% 68754321	. 21
% 68574321	. 21 3
% 68574321	. 21 3 . 
% 36857421	. 21 3 . 54321
% 36857421	. 21 3 . 54321 .
% 13685742	. 21 3 . 54321 . 7654321
% 
% 13685742
% \end{verbatim}
$$
\begin{tinyrc}{
\begin{array}{@{}|@{\,}c@{\,}|@{\,}c@{\,}|@{\,}c@{\,}|@{\,}c@{\,}
		 |@{\,}c@{\,}|@{\,}c@{\,}|@{\,}c@{\,}|@{\,}c@{\,}|@{}}
\hline\cdot&  +  &\cdot&  +  &  +  &\cdot&  +  &
\hln  \cdot&  +  &\cdot&  +  &  +  &\cdot&     &\phantom{+}
\hln  \cdot&  +  &\cdot&  +  &\cdot&     &     &
\hln  \cdot&  +  &\cdot&  +  &     &     &     &
\hln  \cdot&  +  &\cdot&     &     &     &     &
\hln  \cdot&  +  &     &     &     &     &     &
\hln  \cdot&     &\phantom{+}&&    &\phantom{+}&&
\hln \phantom{+}&&     &     &     &     &\phantom{+}&
\\\hline\multicolumn{8}{c}{}
\\\multicolumn{8}{c}{\hbox{\normalsize{$v$}}}
\end{array}
}\end{tinyrc}
\begin{array}{@{\qquad}c@{\qquad}}
  \\\hbox{and}\\\\\mbox{}
\end{array}
\begin{tinyrc}{
\begin{array}{@{}|@{\,}c@{\,}|@{\,}c@{\,}|@{\,}c@{\,}|@{\,}c@{\,}
		 |@{\,}c@{\,}|@{\,}c@{\,}|@{\,}c@{\,}|@{\,}c@{\,}|@{}}
\hline\cdot&  +  &\cdot&  +  &  +  &\cdot&  +  &
\hln  \cdot&  +  &\cdot&  +  &\cdot&\cdot&     &\phantom{+}
\hln  \cdot&  +  &\cdot&  +  &\cdot&     &     &
\hln  \cdot&  +  &\cdot&  +  &     &     &     &
\hln  \cdot&  +  &\cdot&     &     &     &     &
\hln  \cdot&  +  &     &     &     &     &     &
\hln  \cdot&     &\phantom{+}&&    &\phantom{+}&&
\hln \phantom{+}&&     &     &     &     &\phantom{+}&
\\\hline\multicolumn{8}{c}{}
\\\multicolumn{8}{c}{\hbox{\normalsize{$vs_4$}}}
\end{array}
}\end{tinyrc}
$$
in which the row indices of the dots give the lex first reduced
expressions.%
\end{example}

Whether or not a path from $w_0$ to~$w$ is poptotic, breadth-first search
on the preimage tree (ordering the mitosis offspring as in
Proposition~\ref{prop:offspring}) yields a total order on $\rc(w)$.  It
can be shown that {\em poptotic}\/ total orders by breadth-first search
are linear extensions of the partial order on reduced pipe dreams
determined by chute operations.

Define the simplicial complex $\LL_w$ with vertex set $[n] \times [n]$ to
have as its facets the {\em complements}\/ of the reduced pipe dreams
for~$w$:
\begin{eqnarray*}
  \facets(\LL_w) &=& \big\{\big([n] \times [n]\big) \minus D \mid D \in
  \rc(w)\big\}.
\end{eqnarray*}
This is an example of a `subword complex' \cite{grobGeom,KMsubword}, and
is hence shellable by \cite[Theorem~A.4]{grobGeom}.  Through heuristic
arguments and computer calculations in small symmetric groups, we are
convinced of the following.

\begin{conj} \label{conj:poptotic}
Poptotic orders on $\rc(w)$ by breadth-first~search yield
shellings~of~$\LL_w$.
\end{conj}

To emphasize: shellability is not in question, because shellings of
$\LL_w$ appear in \cite[Theorem~A.4 and Section~3.9]{grobGeom}.  The
conjecture would just give more intuitive shellings than those known.  It
is conceivable that all of the apoptotic total orders are shellings, too,
although this seems less likely.

%%%%%%%%%%%%%%%%%%%%%%%%%%%%%%%%%%%%%%%%%%%%%%%%%%%%%%%%%%%%%%%%%%%%%%%%%
%%%%%%%%%%%%%%%%%%%%%%%%%%%%%%%%%%%%%%%%%%%%%%%%%%%%%%%%%%%%%%%%%%%%%%%%%

%\addcontentsline{toc}{subsection}{\numberline{}Acknowledgements}
\bigskip
\noindent
\textbf{Acknowledgements.}  The author is grateful to Allen Knutson,
Cristian Lenart, Vic Reiner, and Anne Schilling for motivation and
inspiration.  Mitosis was discovered in part because of a superb AMS
meeting on Modern Schubert Calculus organized by Frank Sottile and
Nantel Bergeron (the latter also supplied \LaTeX\ macros for drawing
pipe dreams).

%%%%%%%%%%%%%%%%%%%%%%%%%%%%%%%%%%%%%%%%%%%%%%%%%%%%%%%%%%%%%%%%%%%%%%%%%
%%%%%%%%%%%%%%%%%%%%%%%%%%%%%%%%%%%%%%%%%%%%%%%%%%%%%%%%%%%%%%%%%%%%%%%%%
%\bibliographystyle{amsalpha}\bibliography{biblio}\end{document}
\def\cprime{$'$}
\providecommand{\bysame}{\leavevmode\hbox to3em{\hrulefill}\thinspace}

%%%%%%%%%%%%%%%%%%%%%%%%%%%%%%%%%%%%%%%%%%%%%%%%%%%%%%%%%%%%%%%%%%%%%%%%
%%%%%%%%%%%%%%%%%%%%%%%%%%%%%%%%%%%%%%%%%%%%%%%%%%%%%%%%%%%%%%%%%%%%%%%%

%%%%%%%%%%%%%%%%%%%%%%%%%%%%%%%%%%%%%%%%%%%%%%%%%%%%%%%%%%%%%%%%%%%%%%%%%
\end{document}